\documentclass[11pt]{article}
\usepackage{amsfonts,latexsym,rawfonts,amsmath,amssymb,amsthm,graphicx}
\numberwithin{equation}{section}
\newtheorem{theorem}{Theorem}[section]
\newtheorem{lem}[theorem]{Lemma}
\newtheorem{thm}[theorem]{Theorem}
\newtheorem{pro}[theorem]{Proposition}
\newtheorem{cor}[theorem]{Corollary}

\newtheorem{rem}[theorem]{Remark}
\def\s{\,\,\,\,}

\def\endproof{$\hfill\Box$\\}

\title{ \bf Energy identity and removable singularities of maps from a Riemann surface with tension field unbounded in $L^2$}
\author{ Yong Luo\footnote{ The author is supported by the DFG Collaborative Research Center SFB/Transregio 71.}}

\date{}

\begin{document}
\maketitle
\begin{abstract}
We prove the removal singularity results for maps with bounded energy from the unit disk $B$ of $R^2$ centered at the origin to a closed Riemannian manifold whose tension field is unbounded in $L^2(B)$ but satisfies the following condition:
\begin{eqnarray*}
 (\int_{B_t\setminus B_{\frac{t}{2}}}|\tau(u)|^2)^\frac{1}{2}\leq C_1(\frac{1}{t})^a,
\end{eqnarray*}
for some $0<a<1$ and for any $0<t<1$, where $C_1$ is a constant independent of $t$.

We will also prove that if a sequence $\{u_n\}$ has uniformly bounded energy and satisfies
\begin{eqnarray*}
 (\int_{B_t\setminus B_{\frac{t}{2}}}|\tau(u_n)|^2)^\frac{1}{2}\leq C_2(\frac{1}{t})^a,
 \end{eqnarray*}
for some $0<a<1$ and for any $0<t<1$, where $C_2$ is a constant independent of $n$ and $t$, then the energy identity holds for this sequence and there will be no neck formation during the blow up process.
\end{abstract}
\section{Introduction}
  Let $(M,g)$ be a Riemannian manifold and $(N,h)$ be a Riemannian manifold without boundary. For a $W^{1,2}(M,N)$ map $u$, the energy density of $u$ is defined by
$$e(u)=\frac{1}{2}|\nabla u|^2=Tr_g(u^\ast h).$$
where $u^\ast h$ is the pullback of the metric tensor $h$.

The energy functional of the mapping $u$ is defined as
$$E(u)=\int_Me(u)dV.$$

A map $u\in C^1(M,N)$ is called a harmonic map if it is a critical point of the energy.

By Nash embedding theorem $N$ can be isometrically embedded into an Euclidean space $R^K$ for some positive integer $K$. Then $(N,h)$ can be viewed as a submanifold of $R^K$ and a map $u\in W^{1,2}(M,N)$ is a map in $ W^{1,2}(M,R^K)$ whose image lies on $N$. The space $C^1(M,N)$ should be understood in the same way. In this sense we have the following Euler-Lagrangian equation for harmonic maps
$$\triangle u=A(u)(\nabla u,\nabla u).$$

The tension field of a map $u$, $\tau(u)$, is defined by
$$\tau(u)=\triangle u-A(u)(\nabla u,\nabla u),$$
where $A$ is the second fundamental form of $N$ in $R^K$. So $u$ is a harmonic map if and only if $\tau(u)=0$.

Notice that when $M$ is a Riemann surface the functional $E(u)$ is conformal invariant and harmonic maps are of special interest in this case. Consider a harmonic map $u$ from a Riemann surface $M$ to $N$, recall that in the fundamental paper {\cite{S-U}} of Sacks and Uhlenbeck the well-known removable singularity theorem has long been established by using a class of piecewise smooth harmonic functions to approximate the weak harmonic map. In the paper {\cite{Li-W1}}, the authors gave a slightly different proof of the following removable singularity theorem:
\begin{thm}{\cite{Li-W1}}
Let $B$ be the unit disk in $R^2$ centered at the origin. If $u: B\setminus\{0\}\rightarrow N$ is a $W^{2,2}_{loc}(B\setminus\{0\},N)\cap W^{1,2}(B,N)$ map and u satisfies
$$\tau(u)=g\in L^2(B,R^K),$$
then $u$ may be extended to a map belonging to $W^{2,2}(B,N).$
\end{thm}
In this direction we will prove the following result:
\begin{pro}\label{Gradient estimates}
Let $B$ be the unit disk in $R^2$ centered at the origin. If
$$u: B\setminus\{0\}\rightarrow N$$
is a $W^{2,2}_{loc}(B\setminus\{0\},N)\cap W^{1,2}(B,N)$ map and u satisfies the following condition:
$$ (\int_{B_t\setminus B_{\frac{t}{2}}}|\tau(u)|^2)^\frac{1}{2}\leq C(\frac{1}{t})^a,$$
for some $0<a<1$ and for any $0<t<1$, where $C$ is a constant independent of $t$, then there exists some $s>1$ such that
$$\nabla u\in L^{2s}(B).$$
\end{pro}
A direct corollary of this result is the following removable singularity theorem:
\begin{thm}\label{removable-singularity}
Assume that $u\in W^{2,2}_{loc}(B\setminus\{0\},N)\cap W^{1,2}(B,N)$ and u satisfies the following condition:
$$ (\int_{B_t\setminus B_{\frac{t}{2}}}|\tau(u)|^2)^\frac{1}{2}\leq C(\frac{1}{t})^a,$$
for some $0<a<1$ and for any $0<t<1$, where $C$ is a constant independent of $t$.

Then we have
$$u\in\bigcap_{1<p<\frac{2}{1+a}}W^{2,p}(B,N).$$
\end{thm}

Consider a sequence of maps $\{u_n\}$ from a Riemann surface $M$ to $N$ with uniformly bounded energy. It is clear that $\{u_n\}$ converges to $u$ weakly in $W^{1,2}(M,N)$ for some $u\in W^{1,2}(M,N)$, but in general it may not converge strongly in $W^{1,2}(M,N)$ to $u$, and the falling of the strong convergence is due to the energy concentration at finite points. When $\tau(u_n)=0,$ i.e. $u_n$ are harmonic maps, Jost (\cite{J}) and Parker ({\cite{P}) independently proved that the lost energy is exactly the the sum of the energy of the bubbles, recall that in the fundamental paper {\cite{S-U}} Sacks and Uhlenbeck proved the bubbles for such a sequence are harmonic spheres defined as nontrivial harmonic maps from $S^2$ to $N$. This result is called energy identity. Furthermore Parker ({\cite{P}) proved that there is no neck formation during the blow up process, i.e. there holds true the bubble tree convergence.

For the case when $\tau(u_n)$ is bounded in $L^2$, i.e. $\{u_n\}$ is an approximated harmonic map sequence, the energy identity was proved in {\cite{Q}} for N is a sphere and proved in {\cite{D-T}} for the general target manifold $N$(the same result also proved by Wang \cite{W} independently). In {\cite{Q-T}} the authors proved that there is no neck formation during the blow up process. See also \cite{Lin-W1}. For the heat flow of harmonic maps, related results can also be found in {\cite{Top1}}{\cite{Top2}}. If the target manifold is a sphere, the energy identity and bubble tree convergence were proved by Fanghua Lin and Changyou Wang (\cite{Lin-W2}) for sequences with tension fields uniform  bounded in $L^p$, for any $p>1$. In fact, they proved this result under a scaling invariant condition which can be deduced from the uniform boundness of the tension field in  $L^p$.

By virtue of Fanghua Lin and Changyou Wang's result, it is natural to ask the following question:

\textbf{Qestion:} Let $\{u_n\}$ be a sequence from a closed Riemann surface to a closed Riemannian manifold with tension field uniformly bounded in $L^p$ for some $p>1$. Do the energy identity and bubble tree convergence results hold true during blowing up for such a sequence?
\begin{rem}
Parker ({\cite{P}} constructed a sequence from a Riemann surface whose tension field is uniformly bounded in $L^1$, in which the energy identity fails.
\end{rem}
In {\cite{Li-Z}} the authors proved the following theorem:
\begin{thm}{\cite{Li-Z}}\label{Li-Zhu}
Let $\{u_n\}$ be a sequence of maps from $B$ to $N$ in $W^{1,2}(B,N)$ with tension field $\tau(u_n)$, where $B$ is the unit disk of $R^2$ centered at the origin. If
\\(I) $\|u_n\|_{W^{1,2}(B)}$+$\|\tau(u_n)\|_{W^{1,p}(B)}\leq\Lambda$ for some $p\geq\frac{6}{5}$,
\\(II) $u_n\rightarrow u$ strongly in $W^{1,2}_{loc}(B\setminus \{0\},N)$ as $n\rightarrow\infty$,
\\then there exists a subsequence of $\{u_n\}$(still denoted by $\{u_n\}$) and some nonnegative integer $k$, such that for any $i=1,...,k$, there are some points $x_n^i$, positive numbers $r_n^i$ and a nonconstant harmonic sphere $\omega^i$(which is viewed as a map from $R^2\cup\{\infty\}\rightarrow N$) such that
\\(1) $x_n^i\rightarrow 0$, $r_n^i\rightarrow 0$ as $n\rightarrow\infty$;
\\(2) $\lim_{n\rightarrow\infty}(\frac{r_n^i}{r_n^j}+\frac{r_n^j}{r_n^i}+\frac{|x_n^i-x_n^j|}{r_n^i+r_n^j})=\infty$ for any $i\neq j$;
\\(3) $\omega^i$ is the weak limit or strong limit of $u_n(x_n^i+r_n^ix)$ in $W^{1,2}_{loc}(R^2,N)$;
\\(4) Energy identity:
$$\lim_{n\rightarrow\infty}E(u_n,B)=E(u,B)+\sum_{i=1}^kE(\omega^i,R^2).$$
\\(5) Neckless: the image $u(B)\bigcup_{i=1}^k\omega^i(R^2)$ is a connected set.
\end{thm}

\begin{lem}
If $\tau(u)$ satisfies
$$(\int_{B_t\setminus B_{\frac{t}{2}}}|\tau(u)|^2)^\frac{1}{2}\leq C(\frac{1}{t})^a,$$
for some $0<a<1$ and for any $0<t<1$, where $C$ is a constant independent of $t$, then $\tau(u)$ is bounded in $L^p(B)$ for some $p\geq \frac{6}{5}.$
\proof
\begin{eqnarray*}
\int_{B_{2^{-k+1}}\setminus B_{2^{-k}}}|\tau(u)|^p&\leq& C(2^{-k})^{2-p}\|\tau(u)\|^p_{L^2(B_{2^{-k+1}}\setminus B_{2^{-k}})}\\&\leq& C(2^{-k})^{2-p-ap},
 \end{eqnarray*}
 hence $$\int_B|\tau(u)|^p\leq C\sum_{k=1}^\infty(2^{-k})^{2-p-ap}.$$
 When $0<a<\frac{2}{3}$, we can choose some $p\geq\frac{6}{5}$ such that $2-p-ap>0$ and so $\sum_{k=1}^\infty(2^{-k})^{2-p-ap}\leq C$, which implies that $\tau(u)$ is bounded in $L^p(B)$ for some $p\geq\frac{6}{5}$.
\end{lem}
 Thus theorem \ref{Li-Zhu} holds true for sequences $\{u_n\}$ satisfying the following conditions
\\(I) $\|u_n\|_{W^{1,2}(B)}\leq\Lambda$ and $(\int_{B_t\setminus B_{\frac{t}{2}}}|\tau(u_n)|^2)^\frac{1}{2}\leq C(\frac{1}{t})^a$, for some $0<a<\frac{2}{3}$ and for any $0<t<1$, where $C$ is independent of $n$, $t$;
\\(II)  $u_n\rightarrow u$ strongly in $W^{1,2}_{loc}(B\setminus \{0\},N)$ as $n\rightarrow\infty$.

With the help of this observation, we find the following theorem:
\begin{thm}\label{bubble tree}
Let $\{u_n\}$ be a sequence of maps from $B$ to $N$ in $W^{1,2}(B,N)$ with tension field $\tau(u_n)$, where $B$ is the unit disk of $R^2$ centered at the origin. If
\\(I) $\|u_n\|_{W^{1,2}(B)}\leq\Lambda$ and
$$(\int_{B_t\setminus B_{\frac{t}{2}}}|\tau(u_n)|^2)^\frac{1}{2}\leq C(\frac{1}{t})^a,$$
for some $0<a<1$ and for any $0<t<1$, where $C$ is independent of $n$, $t$,
\\(II) $u_n\rightarrow u$ strongly in $W^{1,2}_{loc}(B\setminus \{0\},N)$ as $n\rightarrow\infty$,
\\then there exists a subsequence of $\{u_n\}$(still denote it by $\{u_n\}$) and some nonnegative integer $k$, such that for any $i=1,...,k$, there are some points $x_n^i$, positive numbers $r_n^i$ and a nonconstant harmonic sphere $\omega^i$(which is viewed as a map from $R^2\cup\{\infty\}\rightarrow N$) such that
\\(1) $x_n^i\rightarrow 0$, $r_n^i\rightarrow 0$ as $n\rightarrow\infty$;
\\(2) $\lim_{n\rightarrow\infty}(\frac{r_n^i}{r_n^j}+\frac{r_n^j}{r_n^i}+\frac{|x_n^i-x_n^j|}{r_n^i+r_n^j})=\infty$ for any $i\neq j$;
\\(3) $\omega^i$ is the weak limit or strong limit of $u_n(x_n^i+r_n^ix)$ in $W^{1,2}_{loc}(R^2,N)$;
\\(4) Energy identity:
$$\lim_{n\rightarrow\infty}E(u_n,B)=E(u,B)+\sum_{i=1}^kE(\omega^i,R^2).$$
\\(5) Neckless: the image $u(B)\bigcup_{i=1}^k\omega^i(R^2)$ is a connected set.
\end{thm}
\begin{rem}

When
$$(\int_{B_t\setminus B_{\frac{t}{2}}}|\tau(u_n)|^2)^\frac{1}{2}\leq C(\frac{1}{t})^a,$$
for some $0<a<1$ and for any $0<t<1$, where $C$ is independent of $n$, $t$, we can deduce that $\tau(u_n)$ is uniformly bounded in $L^p(B)$ for any $p<\frac{2}{1+a}$ and when $a\rightarrow1$, $p\rightarrow1$. Hence our condition is stronger than the condition that the tension field is bounded in  $L^p$ , for some  $p>1$, and this result suggests that we probably have a positive answer to the question on page 3.
\end{rem}

This paper is organized as follows: In section 2 we will quote and prove results which will be important in the following sections. In section 3 we will prove the removable singularity result and the last theorem, theorem \ref{bubble tree}, will be proved in section 4.

Throughout this paper, without illustration the letter $C$ will denote positive constants varying from line to line, and we do not always distinguish sequences and its subsequences.

\textbf{Acknowledgement} I would like to thank Professor Guofang Wang for pointing out many typing errors and I also appreciate Professor Youde Wang for bringing my attention to the paper {\cite{Li-W1}}. My interest in this kind of problem began at a class given by Professor Yuxiang Li at Tinghua university, and I had many useful discussions with him .
\section{The $\epsilon$-regularity lemma and the Pohozeav inequality}
In this section we give the well known small energy regularity lemma for approximated harmonic maps and a version of  Pohozeav inequality, which will be important in the following sections.  We assume that the disk $B\subseteq R^2$ is the unit disk centered at the origin, which has the standard flat metric.
\begin{lem}
Suppose that $u\in W^{2,2}(B,N)$ and $\tau(u)=g\in L^2(B,R^K),$ then there exists an $\varepsilon_0>0$ such that if $\int_B|\nabla u|^2\leq \varepsilon_0^2,$  we have
\begin{eqnarray}
\|u-\bar{u}\|_{W^{2,2}(B_{\frac{1}{2}})}\leq C(\|\nabla u\|_{L^2(B)}+\|g\|_{L^2(B)}).
\end{eqnarray}
Here $\bar{u}$ is the mean value of $u$ over $B_{\frac{1}{2}}$.
\end{lem}
\proof We can find a complete proof of this lemma in {\cite{D-T}}.
\endproof

Using the standard elliptic estimates and the embedding theorems we can derive from the above lemma that
\begin{cor}\label{cor1}
Under the assumptions of proposition \ref{Gradient estimates}, we have
\begin{eqnarray}
Osc_{B_{2r}\setminus B_r}u\leq C(\|\nabla u\|_{L^2(B_{4r}\setminus B_{\frac{r}{2}})}+r\|g\|_{L^2(B_{4r}\setminus B_{\frac{r}{2}})})\leq C(\|\nabla u\|_{L^2(B_{4r}\setminus B_{\frac{r}{2}})}+r^{1-a}).
\end{eqnarray}
\end{cor}
In order to prove proposition \ref{Gradient estimates}, we need the following Pohozeav inequality
\begin{lem}
Under the assumptions of proposition \ref{Gradient estimates}, for $0<t_2<t_1<1$, there holds true that
\begin{eqnarray}
\int_{\partial(B_{t_1}\setminus B_{t_2})}r(|\frac{\partial u}{\partial r}|^2-\frac{1}{2}|\nabla u|^2)\leq t_1\|\nabla u\|_{L^2(B_{t_1}\setminus B_{t_2})}\|g\|_{L^2(B_{t_1}\setminus B_{t_2})}.
\end{eqnarray}
\end{lem}
\proof Multiplying the both sides of the equation $\tau(u)=g$ by $r\frac{\partial u}{\partial r}$, we get
\begin{equation*}
\int_{B_{t_1}\setminus B_{t_2}}r\frac{\partial u}{\partial r}\bigtriangleup u=\int_{B_{t_1}\setminus B_{t_2}}r\frac{\partial u}{\partial r}g.
\end{equation*}
By integral by parts we get
$$\int_{B_{t_1}\setminus B_{t_2}}r\frac{\partial u}{\partial r}\bigtriangleup udx=\int_{\partial(B_{t_1}\setminus B_{t_2})}r|\frac{\partial u}{\partial r}|^2-\int_{B_{t_1}\setminus B_{t_2}}\nabla(r\frac{\partial u}{\partial r})\nabla udx,$$
and
\begin{eqnarray*}
\int_{B_{t_1}\setminus B_{t_2}}\nabla(r\frac{\partial u}{\partial r})\nabla u dx&=&\int_{B_{t_1}\setminus B_{t_2}}\nabla(x^k\frac{\partial u}{\partial x^k})\nabla udx\\&=&\int_{B_{t_1}\setminus B_{t_2}}|\nabla u|^2+\int_{t_2}^{t_1}\int_0^{2\pi}\frac{r}{2}\frac{\partial}{\partial r}|\nabla u|^2rd\theta dr\\&=&\int_{B_{t_1}\setminus B_{t_2}}|\nabla u|^2+\frac{1}{2}\int_{\partial(B_{t_1}\setminus B_{t_2})}|\nabla u|^2r-\int_{B_{t_1}\setminus B_{t_2}}|\nabla u|^2\\&=&\frac{1}{2}\int_{\partial(B_{t_1}\setminus B_{t_2})}|\nabla u|^2r.
\end{eqnarray*}
Clearly it implies the conclusion of the lemma.
\endproof

By the above lemma, we can deduce that
\begin{cor}\label{cor2}
Under the assumptions of proposition \ref{Gradient estimates}, we have
$$\frac{\partial}{\partial t}\int_{B_t\setminus B_{\frac{t}{2}}}|\frac{\partial u}{\partial r}|^2-\frac{1}{2}|\nabla u|^2\leq C\|\nabla u\|_{L^2(B_t\setminus B_{\frac{t}{2}})}t^{-a}.$$
\end{cor}
\proof In the previous lemma, let $t_1=t$, $t_2=\frac{t}{2}$, then there holds true
\begin{eqnarray*}
\frac{\partial}{\partial t}\int_{B_t\setminus B_{\frac{t}{2}}}|\frac{\partial u}{\partial r}|^2-\frac{1}{2}|\nabla u|^2&=&
\int_{\partial{B_t}}(|\frac{\partial u}{\partial r}|^2-\frac{1}{2}|\nabla u|^2)-\frac{1}{2}\int_{\partial B_{\frac{t}{2}}}(|\frac{\partial u}{\partial r}|^2-\frac{1}{2}|\nabla u|^2)
\\&\leq&\|g\|_{L^2(B_{t}\setminus B_{\frac{t}{2}})}\|\nabla u\|_{L^2(B_{t}\setminus B_{\frac{t}{2}})}
\\&\leq&C\|\nabla u\|_{L^2(B_{t}\setminus B_{\frac{t}{2}})}t^{-a}.
\end{eqnarray*}
\endproof

A direct corollary of the above conclusion is the following
\begin{cor}
Under the assumptions of proposition \ref{Gradient estimates}, we can get that
\begin{eqnarray}
\int_{B_t\setminus B_{\frac{t}{2}}}|\frac{\partial u}{\partial r}|^2-\frac{1}{2}|\nabla u|^2\leq C\|\nabla u\|_{L^2(B_t)}t^{1-a}.\label{inequ1}
\end{eqnarray}
\end{cor}
\proof Integrating from 0 to $t$ over the both sides of the inequality of the above corollary, and noting that $\|\nabla u\|_{L^2(B_{s}\setminus B_{\frac{s}{2}})}\leq\|\nabla u\|_{L^2(B_t)}$, for any $s\leq t$, we get this inequality directly.
\endproof
\section{Removable of singularities}
In this section we will discuss the issues on the removable of singularities of a class of approximated harmonic maps from the unit disk of $R^2$ centered at the origin to a closed Riemannian manifold $N$.
\begin{lem}
Assume that $u$ satisfies the assumptions of proposition \ref{Gradient estimates}, then there are constants $\lambda>0$, $C>0$, such that the following holds true
\begin{eqnarray}
\int_{B_r}|\nabla u|^2\leq Cr^\lambda,
\end{eqnarray}
for $r$ small enough.
\end{lem}
\proof Because we only need to prove the lemma for $r$ is small, we can assume that $E(u,B)<\varepsilon_0$. Let $u^\ast(r): (0,1)\rightarrow R^K$ be a curve defined as follows
$$u^\ast(r)=\frac{1}{2\pi}\int_0^{2\pi}u(r,\theta)d\theta.$$
Then
$$\frac{\partial u^\ast}{\partial r}=\frac{1}{2\pi}\int_0^{2\pi}\frac{\partial u}{\partial r}d\theta.$$
On the one hand, we have
\begin{eqnarray*}
\int_{B_{2^{-k}t}\setminus B_{2^{-k-1}t}}\nabla u\nabla(u-u^\ast)&\geq&\int_{B_{2^{-k}t}\setminus B_{2^{-k-1}t}}(|\nabla u|^2-|\frac{\partial u}{\partial r}|^2)\\&\geq&\frac{1}{2}\int_{B_{2^{-k}t}\setminus B_{2^{-k-1}t}}|\nabla u|^2-C(2^{-k}t)^{1-a},
\end{eqnarray*}
where in the above second inequality we have used (\ref{inequ1}).

Summing $k$ from 0 to infinity we get
$$\int_{B_t}\nabla u\nabla(u-u^\ast)\geq\frac{1}{2}\int_{B_t}|\nabla u|^2-Ct^{1-a}.$$
On the other hand, we have
\begin{eqnarray*}
\int_{B_{2^{-k}t}\setminus B_{2^{-k-1}t}}\nabla u\nabla(u-u^\ast)&=&-\int_{B_{2^{-k}t}\setminus B_{2^{-k-1}t}}(u-u^\ast)\triangle u\\&+&\int_{\partial(B_{2^{-k}t}\setminus B_{2^{-k-1}t})}\frac{\partial u}{\partial r}(u-u^\ast)\\&=&-\int_{B_{2^{-k}t}\setminus B_{2^{-k-1}t}}(u-u^\ast)(\tau(u)-A(u)(\nabla u,\nabla u))\\&+&\int_{\partial(B_{2^{-k}t}\setminus B_{2^{-k-1}t})}\frac{\partial u}{\partial r}(u-u^\ast).
\end{eqnarray*}
Hence by summing $k$ from 0 to infinity we get
\begin{eqnarray*}
\int_{B_t}\nabla u\nabla(u-u^\ast)&\leq&\sum_{k=0}^\infty\|u-u^\ast\|_{L^\infty(B_{2^{-k}t}\setminus B_{2^{-k-1}t})}[\|A\|_{L^\infty}\int_{B_{2^{-k}t}\setminus B_{2^{-k-1}t}}|\nabla u|^2\\&+&C(2^{-k}t)^{1-a}]+\int_{\partial B_t}\frac{\partial u}{\partial r}(u-u^\ast)\\&\leq&\epsilon\int_{B_t}|\nabla u|^2+Ct^{1-a}+\int_{\partial B_t}\frac{\partial u}{\partial r}(u-u^\ast),
\end{eqnarray*}
note that here we have used corollary \ref{cor1} and we let $\epsilon$ small by letting $t$ small.

Note that
\begin{eqnarray*}
|\int_{\partial B_t}\frac{\partial u}{\partial r}(u-u^\ast)|&\leq&(\int_{\partial B_t}|\frac{\partial u}{\partial r}|^2)^\frac{1}{2}(\int_{\partial B_t}|u-u^\ast|^2)^\frac{1}{2}\\&\leq&(\int_0^{2\pi}t^2|\frac{\partial u}{\partial r}|^2d\theta)^\frac{1}{2}(\int_0^{2\pi}|\frac{\partial u}{\partial \theta}|^2d\theta)^\frac{1}{2}\\&\leq&\frac{1}{2}\int_0^{2\pi}(|\frac{\partial u}{\partial \theta}|^2+t^2|\frac{\partial u}{\partial r}|^2)d\theta\\&=&\frac{t}{2}\int_{\partial B_t}|\nabla u|^2.
\end{eqnarray*}
Combining the two sides of inequalities and let $\epsilon$ small(we can do this by letting $t$ small), we conclude that there is a $\lambda$ which is a positive constant smaller than one such that
$$\lambda\int_{B_t}|\nabla u|^2\leq t\int_{\partial B_t}|\nabla u|^2+Ct^{1-a}.$$
Now denote $f(t)=\int_{B_t}|\nabla u|^2$, then we get the following ordinary differential inequality
$$(\frac{f(t)}{t^\lambda})'\geq-Ct^{-\lambda-a}.$$
Note that we can let $\lambda$ small enough such that $\lambda+a<1$, thus we can get that
$$f(t)=\int_{B_t}|\nabla u|^2\leq Ct^\lambda,\s for\s t\s small\s enough.$$\endproof

Now we are in a position to give a complete proof of proposition \ref{Gradient estimates}. Let $r_k=2^{-k}$ and $v_k(x)=u(r_kx)$, then
\begin{eqnarray*}
(\int_{B_2\setminus B_1}|\nabla v_k|^{2s})^\frac{1}{2s}&\leq&C\|v_k-\bar{v_k}\|_{W^{2,2}(B_2\setminus B_1)}\\&\leq&(\int_{B_4\setminus B_{\frac{1}{2}}}|\nabla v_k|^2)^\frac{1}{2}+C(\int_{B_{4r_k}\setminus B_{\frac{1}{2}r_k}}r_k^2|\tau|^2)^\frac{1}{2}.
\end{eqnarray*}
Therefore we deduce that
\begin{eqnarray*}
\int_{B_2\setminus B_1}|\nabla v_k|^{2s}&\leq& C(\int_{B_4\setminus B_{\frac{1}{2}}}|\nabla v_k|^2)^s+C(\int_{B_{4r_k}\setminus B_{\frac{1}{2}r_k}}r_k^2|\tau|^2)^s\\&\leq&C(\int_{B_4\setminus B_{\frac{1}{2}}}|\nabla v_k|^2)^s+Cr_k^{2s(1-a)}.
\end{eqnarray*}
Note that when $k$ is large enough, there holds true
$$\int_{B_{4r_k}\setminus B_{\frac{r_k}{2}}}|\nabla u|^2\leq1.$$
Hence we have that
\begin{eqnarray*}
r_k^{2s-2}\int_{B_{2r_k}\setminus B_{r_k}}|\nabla u|^{2s}&\leq&C(\int_{B_{4r_k}\setminus B_{\frac{r_k}{2}}}|\nabla u|^2)^s+Cr_k^{2s(1-a)}\\&\leq&C\int_{B_{4r_k}\setminus B_{\frac{r_k}{2}}}|\nabla u|^2+Cr_k^{2s(1-a)}.
\end{eqnarray*}
This implies that
$$\int_{B_{2r_k}\setminus B_{r_k}}|\nabla u|^{2s}\leq Cr_k^{2-2s}r_k^\lambda+Cr_k^{2-2sa}.$$
Now choose $s>1$ such that $2s-2<\frac{\lambda}{2}$ and $2-2sa>0$. There exists a positive integer $k_0$ such that when $k\geq k_0$ there holds true
$$\int_{B_{2^{-k+1}}\setminus B_{2^{-k}}}|\nabla u|^{2s}\leq C(2^{\frac{-\lambda}{2}k}+2^{-k(2-2sa)}).$$
Therefore
$$\int_{B_r}|\nabla u|^{2s}\leq C\sum_{k=k_0}^\infty(2^{\frac{-\lambda}{2}k}+2^{-k(2-2sa)}))\leq C,$$
for any $r\leq 2^{-k_0+1}$, and this completes the proof of proposition \ref{Gradient estimates}.

Proof of theorem \ref{removable-singularity}: Note that
\begin{eqnarray*}
\int_{B_{2^{-k}}\setminus B_{2^{-k-1}}}|\tau(u)|^p&\leq&C(2^{-k})^{2-p}(\int_{B_{2^{-k}}\setminus B_{2^{-k-1}}}|\tau(u)|^2)^\frac{p}{2}\\&\leq&C(2^{-k})^{2-p-pa}.
\end{eqnarray*}
Hence we can deduce from this by summing $k$ from 0 to infinity that
$$\int_B|\tau(u)|^p\leq C\s for\s p\s<\s \frac{2}{1+a}.$$
Recall that we have proved that $\nabla u\in L^{2s}(B)$ for some $s>1$, hence we can deduce that
$$u\in\bigcap_{1<p<\frac{2}{1+a}}W^{2,p}(B,N),$$
by standard elliptic estimates and the bootstrap argument.
\endproof
\section{The bubble tree structure}
\subsection{Energy identity}
Assume that $\{u_n\}$ is a sequence which is uniformly bounded in $W^{1,2}(B,N)$ and satisfies
$$(\int_{B_t\setminus B_{\frac{t}{2}}}|\tau(u_n)|^2)^\frac{1}{2}\leq C(\frac{1}{t})^a,$$
for some $0<a<1$ and for any $0<t<1$, where $C$ is a constant independent of $n$ and $t$.

In this section, we will prove the energy identity for this sequence and for convenience we will assume that there is only one bubble $\omega$ which is the strong limit of $u_n(r_n.)$ in $W^{1,2}_{loc}(R^2,N).$

Because
$$\lim_{n\rightarrow\infty}\int_B|\nabla u_n|^2=\lim_{n\rightarrow\infty}\int_{B\setminus B_{\delta}}|\nabla u_n|^2+\lim_{n\rightarrow\infty}\int_{B_\delta\setminus B_{Rr_n}}|\nabla u_n|^2+\lim_{n\rightarrow\infty}\int_{ B_{Rr_n}}|\nabla u_n|^2.$$
and
$$\lim_{\delta\rightarrow0}\lim_{n\rightarrow\infty}\int_{B\setminus B_{\delta}}|\nabla u_n|^2=\int_B|\nabla u|^2,$$
$$\lim_{R\rightarrow\infty}\lim_{n\rightarrow\infty}\int_{ B_{Rr_n}}|\nabla u_n|^2=\int_{R^2}|\nabla\omega|^2,$$
hence to prove the energy identity we only need to prove that
\begin{eqnarray}
\lim_{R\rightarrow\infty}\lim_{\delta\rightarrow0}\lim_{n\rightarrow\infty}\int_{B_\delta\setminus B_{Rr_n}}|\nabla u_n|^2=0.
\end{eqnarray}

At first we note that under the assumption of only one bubble we can deduce by a standard blow up argument the following
\begin{lem}
For any $\epsilon>0$, there exist $R$ and $\delta$ such that
\begin{eqnarray}
\int_{B_{2\lambda}\setminus B_{\lambda}}|\nabla u_n|^2\leq \epsilon^2,\s for\s any\s \lambda\in(\frac{Rr_n}{2},2\delta).
\end{eqnarray}
\end{lem}
Now we have been prepared to prove the energy identity. The proof is a little similar with the proof in the previous section. We assume that $\delta=2^{m_n}Rr_n$ where $m_n$ is a positive integer.

On the one hand, we have
\begin{eqnarray*}
\int_{B_{2^{k}Rr_n}\setminus B_{2^{k-1}Rr_n}}\nabla u_n\nabla(u_n-u_n^\ast)&\geq&\int_{B_{2^{k}Rr_n}\setminus B_{2^{k-1}Rr_n}}(|\nabla u_n|^2-|\frac{\partial u_ n}{\partial r}|^2)\\&\geq&\frac{1}{2}\int_{B_{2^{k}Rr_n}\setminus B_{2^{k-1}Rr_n}}|\nabla u_n|^2-C(2^{k}Rr_n)^{1-a}.
\end{eqnarray*}
This implies that
$$\int_{B_{\delta}\setminus B_{Rr_n}}\nabla u_n\nabla(u_n-u_n^\ast)\geq\frac{1}{2}\int_{B_{\delta}\setminus B_{Rr_n}}|\nabla u_n|^2-C\delta^{1-a}.$$

On the other hand, we have that
\begin{eqnarray*}
\int_{B_{2^{k}Rr_n}\setminus B_{2^{k-1}Rr_n}}\nabla u_n\nabla(u_n-u_n^\ast)&=&-\int_{B_{2^{k}Rr_n}\setminus B_{2^{k-1}Rr_n}}(u_n-u_n^\ast)\triangle u_n\\&+&\int_{\partial(B_{2^{k}Rr_n}\setminus B_{2^{k-1}Rr_n})}\frac{\partial u_n}{\partial r}(u_n-u_n^\ast)\\&=&-\int_{B_{2^{k}Rr_n}\setminus B_{2^{k-1}Rr_n}}(u_n-u_n^\ast)(\tau(u_n)-A(u_n)(\nabla u_n,\nabla u_n))\\&+&\int_{\partial(B_{2^{k}Rr_n}\setminus B_{2^{k-1}Rr_n})}\frac{\partial u_n}{\partial r}(u_n-u_n^\ast).
\end{eqnarray*}
We deduce from the above by summing from 1 to $m_n$ that
\begin{eqnarray*}
\int_{B_\delta\setminus B_{Rr_n}}\nabla u_n\nabla(u_n-u_n^\ast)&\leq&\sum_{k=1}^{m_n}\|u_n-u_n^\ast\|_{L^\infty(B_{2^{k}Rr_n}\setminus B_{2^{k-1}Rr_n})}[\|A\|_{L^\infty}\int_{B_{2^{k}Rr_n}\setminus B_{2^{k-1}Rr_n}}|\nabla u_n|^2\\&+&C(2^{k}Rr_n)^{1-a}]+\int_{\partial (B_\delta\setminus B_{Rr_n})}\frac{\partial u_n}{\partial r}(u_n-u_n^\ast)\\&\leq&\epsilon\int_{B_\delta\setminus B_{Rr_n}}|\nabla u_n|^2+C\delta^{1-a}+\int_{\partial (B_\delta\setminus B_{Rr_n})}\frac{\partial u_n}{\partial r}(u_n-u_n^\ast).
\end{eqnarray*}

Comparing with the two sides we get that
$$(1-2\epsilon)\int_{B_\delta\setminus B_{Rr_n}}|\nabla u_n|^2\leq C\delta^{1-a}+2\int_{\partial (B_\delta\setminus B_{Rr_n})}\frac{\partial u_n}{\partial r}(u_n-u_n^\ast).$$
As for the boundary terms, we have
\begin{eqnarray*}
\int_{\partial B_\delta}\frac{\partial u_n}{\partial r}(u_n-u_n^\ast)&\leq&(\int_{\partial B_\delta}|\frac{\partial u_n}{\partial r}|^2)^\frac{1}{2}(\int_{\partial B_\delta}|u_n-u_n^\ast|^2)^\frac{1}{2}\\&\leq&(\int_0^{2\pi}\delta^2|\frac{\partial u_n}{\partial r}d\theta|^2)^\frac{1}{2}(\int_0^{2\pi}|\frac{\partial u_n}{\partial\theta}|^2d\theta)^\frac{1}{2}\\&\leq&\frac{1}{2}\int_0^{2\pi}\delta^2|\frac{\partial u_n}{\partial r}d\theta|^2+|\frac{\partial u_n}{\partial\theta}|^2d\theta\\&=&\frac{\delta^2}{2}\int_0^{2\pi}|\nabla u_n|^2d\theta.
\end{eqnarray*}
By trace embedding theorem, we have
\begin{eqnarray*}
\int_0^{2\pi}|\nabla u_n(.,\delta)|^2\delta d\theta&=&\int_{\partial B_\delta}|\nabla u_n(.,\delta)|^2dS_\delta
\\&\leq& C\delta\|\nabla u_n\|_{W^{1,2}(B_\frac{3\delta}{2}\setminus B_{\frac{\delta}{2}})}^2
\\&\leq&C\delta\|u_n-\bar{u}_n\|_{W^{2,2}(B_\frac{3\delta}{2}\setminus B_{\frac{\delta}{2}})}^2
\\&\leq&C\delta(\frac{1}{\delta}\|\nabla u_n\|_{L^2(B_{2\delta})}^2+\|\tau(u_n)\|_{L^2(B_{2\delta}\setminus B_{\frac{\delta}{4}})}^2)
\\&\leq&C\delta^{1-2a},
\end{eqnarray*}
for $\delta$ small. From this we deduce that
 $$\int_{\partial B_\delta}\frac{\partial u_n}{\partial r}(u_n-u_n^\ast)\leq C\delta^{2(1-a)}.$$
  Similarly we get
 $$\int_{\partial B_{Rr_n}}\frac{\partial u_n}{\partial r}(u_n-u_n^\ast)\leq C(Rr_n)^{2(1-a)},$$
 for $n$ big enough. Therefore
  $$(1-2\epsilon)\int_{B_\delta\setminus B_{Rr_n}}|\nabla u_n|^2\leq C\delta^{1-a}+C\delta^{2(1-a)}+C(Rr_n)^{2(1-a)}.$$
 Clearly this implies that
\begin{eqnarray}
\lim_{R\rightarrow\infty}\lim_{\delta\rightarrow0}\lim_{n\rightarrow\infty}\int_{B_\delta\setminus B_{Rr_n}}|\nabla u_n|^2=0,
\end{eqnarray}
which completes the proof of the energy identity.
\endproof
\subsection{Neckless}
In this part we will prove that there is no neck between the base map $u$ and the bubble $\omega$, i.e. the $C^0$ compactness of the sequence modulo bubbles.

We only need to prove the following
\begin{eqnarray}
\lim_{R\rightarrow\infty}\lim_{\delta\rightarrow0}\lim_{n\rightarrow\infty}Osc_{B_\delta\setminus B_{Rr_n}}u_n=0.
\end{eqnarray}

Again we assume that $\delta=2^{m_n}Rr_n$ and let $Q(t)=B_{2^{t+t_0}Rr_n}\setminus B_{2^{t_0-t}Rr_n}$, similarly to the proof of the previous part we can get
\begin{eqnarray*}
(1-2\epsilon)\int_{Q(k)}|\nabla u_n|^2&\leq& 2^{k+t_0}Rr_n\int_{\partial B_{2^{k+t_0}Rr_n}}|\nabla u_n|^2+2^{t_0-k}Rr_n\int_{\partial B_{2^{t_0-k}Rr_n}}|\nabla u_n|^2\\&+&C(2^{k+t_0}Rr_n)^{1-a}.
\end{eqnarray*}

Denote $f(t)=\int_{Q(t)}|\nabla u_n|^2,$ then we have that
\begin{eqnarray*}
(1-2\epsilon)f(t)&\leq& (1-2\epsilon)f(k+1)\\&\leq&\frac{1}{\log2}f'(k+1)+C(2^{k+t_0}Rr_n)^{1-a},
\end{eqnarray*}
for $k\leq t\leq k+1$.

Note that
\begin{eqnarray*}
f'(k+1)-f'(t)&=&\int_{\partial(B_{2^{k+1+t_0}Rr_n}\setminus B_{2^{t+t_0}Rr_n})}\frac{\partial u_n}{\partial r}(u_n-u_n^\ast)+\int_{\partial(B_{2^{t_0-t}Rr_n}\setminus B_{2^{t_0-k-1}Rr_n})}\frac{\partial u_n}{\partial r}(u_n-u_n^\ast)\\&\leq&C(2^{t+t_0}Rr_n)^{2(1-a)}.
\end{eqnarray*}

Then we have
\begin{eqnarray}
(1-2\epsilon)f(t)\leq\frac{1}{\log2}f'(t)+C(2^{t+t_0}Rr_n)^{1-a}.
\end{eqnarray}
In virtue of the above inequality we get
\begin{eqnarray*}
(2^{-(1-2\epsilon)t}f(t))'&=&2^{-(1-2\epsilon)t}f'(t)-(1-2\epsilon)2^{-(1-2\epsilon)t}f(t)\log2\\&\geq&-C2^{(1-a-(1-2\epsilon))t}(2^{t_0}Rr_n)^{1-a}.
\end{eqnarray*}
Integrating this from 1 to $L$ we get
\begin{eqnarray*}
2^{-(1-2\epsilon)L}f(L)-2^{-(1-2\epsilon)}f(1)&\geq&-C\int_1^L2^{(1-a-(1-2\epsilon))t}(2^{t_0}Rr_n)^{1-a}\\&=&-C\frac{2^{(1-a-(1-2\epsilon))t}}{\log2(1-a-(1-2\epsilon))}|_1^L
(2^{t_0}Rr_n)^{1-a}\\&\geq&-C(2^{t_0}Rr_n)^{1-a}.
\end{eqnarray*}
Therefore we have that
\begin{eqnarray}
f(1)\leq f(L)2^{-(1-2\epsilon)(L-1)}+C(2^{t_0}Rr_n)^{1-a}.
\end{eqnarray}

Now let $t_0=i$ and $D_i=B_{2^{i+1}Rr_n}\setminus B_{2^iRr_n}$, then we have
$$f(1)=\int_{D_i\bigcup D_{i-1}}|\nabla u_n|^2,$$
and equality will hold true for $L$ satisfying
$$Q(L)\subseteq B_\delta\setminus B_{Rr_n}=B_{2^{m_n}Rr_n}\setminus B_{Rr_n}.$$
In other words $L$ should satisfy
\\(1) $i-L\geq 0$;
\\(2) $i+L\leq M_n.$

If
\\(I) $i\leq \frac{1}{2}m_n$, let $L=i$ and so $$f(1)=\int_{D_i\bigcup D_{i-1}}|\nabla u_n|^2\leq CE^2(u_n,B_\delta\setminus B_{Rr_n})2^{-(1-2\epsilon)i}+C(2^iRr_n)^{1-a};$$
\\ (II) $i>\frac{1}{2}m_n$, let $L=m_n-i$ and so $$f(1)=\int_{D_i\bigcup D_{i-1}}|\nabla u_n|^2\leq CE^2(u_n,B_\delta\setminus B_{Rr_n})2^{-(1-2\epsilon)(m_n-i)}+C(2^iRr_n)^{1-a}.$$
Hence we have
\begin{eqnarray*}
\sum_{i=1}^{m_n}E(u_n,D_i)&\leq& CE(u_n,B_\delta\setminus B_{Rr_n})(\sum_{i\leq\frac{1}{2}m_n}2^{-i\frac{1-2\epsilon}{2}}+\sum_{i>\frac{1}{2}m_n}2^{-(m_n-i)\frac{1-2\epsilon}{2}})+C\sum_{i=1}^{m_n}(2^iRr_n)^{\frac{1-a}{2}}\\&\leq&
CE(u_n,B_\delta\setminus B_{Rr_n})+C\delta^\frac{1-a}{2}.
\end{eqnarray*}
Thus we get
\begin{eqnarray*}
Osc_{B_\delta\setminus B_{Rr_n}}u_n&\leq&C\sum_{i=1}^{m_n}(E(u_n,D_i)+(2^iRr_n)^{1-a})\\&\leq&CE(u_n,B_\delta\setminus B_{Rr_n})+C\delta^{\frac{1-a}{2}}.
\end{eqnarray*}

Clearly the above inequality implies that
$$\lim_{R\rightarrow\infty}\lim_{\delta\rightarrow0}\lim_{n\rightarrow\infty}Osc_{B_\delta\setminus B_{Rr_n}}u_n=0,$$
which completes the proof of the Neckless.
\endproof

{}

\vspace{1cm}\sc

YONG LUO

Mathematisches Institut, Albert-Ludwigs-Universit

Eckerstr. 1, 79104, Freiburg, Germany

{\tt yong.luo@math.uni-freiburg.de}
\end{document}